\documentclass[11pt]{article}

\usepackage{amscd,amsmath,amssymb,fancyhdr,epsfig,graphicx,wrapfig}
\usepackage[matrix,arrow]{xy}

\numberwithin{equation}{section}

\newcommand{\version}{version 5.0,\ \ October 08, 2012}

\renewcommand{\leftmark}%
{{\scriptsize Moduli of polarized hyperk\"ahler manifolds}}

\def\eqref#1{(\ref{#1})}

\newcommand{\arrow}{{\:\longrightarrow\:}}
\newcommand{\Z}{{\mathbb Z}}

\def\1{\sqrt{-1}\:}

\newcommand{\cntrct}                % contraction with a vector field
{\hspace{2pt}\raisebox{1pt}{\text{$\lrcorner$}}\hspace{2pt}}

\makeatletter
\def\x@arrow{\DOTSB\Relbar}
\def\xlongequalsignfill@{\arrowfill@\x@arrow\Relbar\x@arrow}
\newcommand{\xlongequal}[2][]{%
        \ext@arrow 0099\xlongequalsignfill@{#1}{#2}}
\def\xlongrightarrowfill@{\arrowfill@\relbar\relbar\longrightarrow}
\newcommand{\xlongrightarrow}[2][]{%
        \ext@arrow 0099\xlongrightarrowfill@{#1}{#2}}
\makeatother

\renewcommand{\bar}{\overline}
\renewcommand{\phi}{\varphi}
\renewcommand{\epsilon}{\varepsilon}
\renewcommand{\geq}{\geqslant}
\renewcommand{\leq}{\leqslant}

\newcommand{\St}{\operatorname{St}}

\newcommand{\Aut}{\operatorname{Aut}}

\newcommand{\rk}{\operatorname{rk}}

\renewcommand{\Re}{\operatorname{Re}}
\renewcommand{\Im}{\operatorname{Im}}

\newcommand{\Comp}{\operatorname{Comp}}
\newcommand{\Teich}{\operatorname{Teich}}
\newcommand{\Diff}{\operatorname{Diff}}
\newcommand{\Mod}{\operatorname{Mod}}
\newcommand{\Pic}{\operatorname{Pic}}
\newcommand{\Per}{\operatorname{\sf Per}}
\newcommand{\Perspace}{\operatorname{{\mathbb P}\sf er}}

\newcounter{Mycounter}[section]
\newcounter{lemma}[section]
\renewcommand{\thelemma}{{Lemma \thesection.\arabic{lemma}}}
\newcommand{\lemma}{%
    \setcounter{lemma}{\value{Mycounter}}
    \refstepcounter{lemma}
    \stepcounter{Mycounter}
    {\noindent \bf \thelemma:\ }}

\newcounter{claim}[section]
\renewcommand{\theclaim}{{Claim \thesection.\arabic{claim}}}
\newcommand{\claim}{%
    \setcounter{claim}{\value{Mycounter}}
    \refstepcounter{claim}
    \stepcounter{Mycounter}
    {\noindent \bf \theclaim:\ }}

\newcounter{sublemma}[section]
\newcounter{corollary}[section]
\renewcommand{\thecorollary}{{Corollary \thesection.\arabic{corollary}}}
\newcommand{\corollary}{%
    \setcounter{corollary}{\value{Mycounter}}
    \refstepcounter{corollary}
    \stepcounter{Mycounter}
    {\noindent \bf \thecorollary:\ }}

\newcounter{theorem}[section]
\renewcommand{\thetheorem}{{Theorem \thesection.\arabic{theorem}}}
\newcommand{\theorem}{%
    \setcounter{theorem}{\value{Mycounter}}
    \refstepcounter{theorem}
    \stepcounter{Mycounter}
    {\noindent \bf \thetheorem:\ }}

\newcommand{\alttheorem}{%
    \setcounter{theorem}{\value{Mycounter}}
    \refstepcounter{theorem}
    \stepcounter{Mycounter}
    {\noindent \bf \thetheorem}}

\newcounter{conjecture}[section]
\newcounter{proposition}[section]
\renewcommand{\theproposition}
      {{Proposition \thesection.\arabic{proposition}}}
\newcommand{\proposition}{%
    \setcounter{proposition}{\value{Mycounter}}
    \refstepcounter{proposition}
    \stepcounter{Mycounter}
    {\noindent \bf \theproposition:\ }}

\newcounter{definition}[section]
\renewcommand{\thedefinition}
      {{Definition~\thesection.\arabic{definition}}}
\newcommand{\definition}{%
    \setcounter{definition}{\value{Mycounter}}
    \refstepcounter{definition}
    \stepcounter{Mycounter}
    {\noindent \bf \thedefinition:\ }}

\newcounter{example}[section]
\newcounter{remark}[section]
\renewcommand{\theremark}{{Remark \thesection.\arabic{remark}}}
\newcommand{\remark}{%
    \setcounter{remark}{\value{Mycounter}}
    \refstepcounter{remark}
    \stepcounter{Mycounter}
    {\noindent \bf \theremark:\ }}

\newcounter{problem}[section]
\newcounter{question}[section]
\makeatletter

\@addtoreset{equation}{section} \@addtoreset{footnote}{section} \makeatother

\def\blacksquare{\hbox{\vrule width 5pt height 5pt depth 0pt}}
\def\endproof{\blacksquare}

\def\SO{\mathop{\text{\rm SO}}}
\def\Gr{\mathop{\text{\rm Gr}}}
\def\O{\mathop{\text{\rm O}}}
\def\GL{\mathop{\text{\rm GL}}}
\def\rk{\mathop{\text{\rm rk}}}
\def\PSL{\mathop{\text{\rm PSL}}}

\begin{document}
%%%%%%%%%%%%%%%%%%%%%%%%%%%%%%%%%%%%%%%%%%%%%%%%%%%%%%%%%%%%
\begin{center}
{\Large\bf
Any component of moduli\\[3mm]
of polarized hyperk\"ahler manifolds\\[3mm]
is dense in its deformation space\\[3mm]
}
%%%%%%%%%%%%%%%%%%%%%%%%%%%%%%%%%%%%%%%%%%%%%%%%%%%%%%%%%%%%

Sasha Anan$'$in, Misha Verbitsky\footnote{Partially supported by the
 RFBR grant 10-01-93113-NCNIL-a, Science Foundation of 
the SU-HSE award No. 10-09-0015 and AG Laboratory HSE, 
RF government grant, ag. 11.G34.31.0023}

\end{center}

%%%%%%%%%%%%%%%%%%%%%%%%%%%%%%%%%%%%%%%%%%%%%%%%
{\small\hspace{0.10\linewidth}
\begin{minipage}[t]{0.84\linewidth}
{\bf Abstract}\\
Let $M$ be a compact hyperk\"ahler manifold, and $W$ the coarse moduli
of complex deformations of $M$. Every positive integer class $v$ in
$H^2(M)$ defines a divisor $D_v$ in $W$ consisting of all algebraic
manifolds polarized by $v$. We prove that every connected component of
this divisor is dense in $W$.
\end{minipage}
}
%%%%%%%%%%%%%%%%%%%%%%%%%%%%%%%%%%%%%%%%%%%%%%%%

\tableofcontents

%%%%%%%%%%%%%%%%%%%%%%%%%%%%%%%%%%%%%%%%%%%%%%%%
\section{Introduction}
%%%%%%%%%%%%%%%%%%%%%%%%%%%%%%%%%%%%%%%%%%%%%%%%

%%%%%%%%%%%%%%%%%%%%%%%%%%%%%%%%%%%%%%%%%%%%%%%%%%%%%%%%%%%%%%%
\subsection{Hyperk\"ahler manifolds and moduli spaces}
%%%%%%%%%%%%%%%%%%%%%%%%%%%%%%%%%%%%%%%%%%%%%%%%%%%%%%%%%%%%%%%

Throughout this paper, a {\bf hyperk\"ahler manifold} means a ``compact
complex manifold admitting a K\"ahler structure and a holomorphically
symplectic form.'' A hyperk\"ahler manifold $M$ is called {\bf simple}
if $\pi_1(M)=0$ and $H^{2,0}(M)=\mathbb C$. By Bogomolov's theorem (see
\cite{_Besse:Einst_Manifo_} and \cite{_Bogomolov:decompo_}), any
hyperk\"ahler manifold has a finite covering which is a product of
simple hyperk\"ahler manifolds and compact tori. Throughout this paper,
we shall silently assume that all our hyperk\"ahler manifolds are
simple. The results that we prove can be stated and proven for general
hyperk\"ahler manifolds, but to do so would destroy the clarity of the
exposition.

For a background story on hyperk\"ahler manifolds, their construction,
and properties, please see \cite{_Beauville_} and
\cite{_Besse:Einst_Manifo_}. The moduli spaces of hyperk\"ahler
manifolds are discussed at great length in \cite{_V:Torelli_}.

The moduli space of complex structures on a given smooth oriented
manifold $M$ is defined, following Kodaira and Spencer, as the quotient
of the Fr\`echet manifold of all integrable almost complex structures
$\Comp$ by the action of the group of orientation-preserving
diffeomorphisms $\Diff^+$, which is considered as a Fr\`echet Lie
group. We denote by $\Comp_0\subset\Comp$ the open set consisting of
all complex structures on $M$ admitting a compatible K\"ahler metric
and a compatible holomorphically symplectic structure. The quotient
$\Mod:=\Comp_0/\Diff^+$ is called a {\bf coarse moduli space} of
hyperk\"ahler manifolds. It is a complex analytic space, usually
non-Hausdorff.

It is well known that a generic point $I\in\Mod$ corresponds to a
non-algebraic complex structure on $M$. In fact, the manifold $(M,I)$
has no divisors, because the corresponding Neron-Severi group
$H^{1,1}(M,\mathbb Z):=H^{1,1}(M)\break\cap H^2(M,\mathbb Z)$ is zero
(see \cite{_Fujiki:HK_}). The algebraic points of $\Mod$ sit on a
countable union of divisors in $\Mod$, which is known to be dense in
$\Mod$ (\cite{_Fujiki:HK_}, \cite{_Verb:alge_}).

In this paper we prove that each of these divisors is itself dense in
$\Mod$. This result is known when $M$ is a K3 surface (this follows
from a statement known as ``Eichler Criterion''; see
\ref{_K3_Remark_}).

%%%%%%%%%%%%%%%%%%%%%%%%%%%%%%%%%%%%%%%%%%%%%%%%%%%%%%%%%%%%%%%
\subsection{Lelong numbers, SYZ conjecture and Gromov's 
precompactness theorem}
%%%%%%%%%%%%%%%%%%%%%%%%%%%%%%%%%%%%%%%%%%%%%%%%%%%%%%%%%%%%%%%

The original motivation for this work came form a research on the
so-called hyperk\"ahler SYZ conjecture (\cite{_Verbitsky:SYZ_}). This
conjecture, which is a version of a (more general) abundance conjecture
of Kawamata, states that a nef bundle on a hyperk\"ahler manifold is
semiample. More specifically, one is interested in holomorphic line
bundles $L$ which are nef, and for which the Bogomolov-Beauville-Fujiki
square of $c_1(L)$ vanishes: $q(c_1(L),c_1(L))=0$ (for a definition of
Bogomolov-Beauville-Fujiki form, see Subsection
\ref{_BBF_Subsection_}). Such line bundles are called {\bf parabolic}.
Any nef bundle admits a singular metric with semipositive curvature
(this follows from general results on weak compactness of positive
currents). If this metric is not ``very singular'', $L$ is effective
(\cite{_Verbitsky:SYZ_}, \cite{_Verbitsky:parabolic_}). The ``not very
singular'' above refers to the vanishing of the so-called Lelong
numbers of the curvature current; these numbers, defined for positive
closed $(p,p)$-currents, vanish for all smooth currents, and measure
the geometric complexity of its singularities in the general case,
taking values in $\mathbb R^{\geq 0}$.

The Lelong numbers are known to be upper semicontinuous in the current
topology. This means, in particular, that any cohomology class $\eta$
which is represented as a limit of currents with Lelong numbers bounded
from below would have positive Lelong numbers.

Suppose now that $\eta\in H^{1,1}(M,\mathbb R)$ is a nef class on a
non-algebraic hyperk\"ahler manifold satisfying $q(\eta,\eta)=0$ (such
class is also called {\bf parabolic}). It is proven in
\cite{_Verbitsky:parabolic_} that the Lelong sets (sets where the
Lelong numbers are bounded from below by a positive number) of $\eta$
are coisotropic with respect to the holomorphic symplectic structure.
However, all complex subvarieties of a generic non-algebraic
hyperk\"ahler manifold are
hyperk\"ahler~(\cite{_Verbitsky:Symplectic_II_}), hence they cannot be
coisotropic. This means that any parabolic nef current on a generic
non-algebraic manifold has vanishing Lelong numbers.

To apply this argument, we need to approximate a given non-algebraic
manifold with a nef current by a sequence of algebraic manifolds with a
rational parabolic current, in a controlled way. To keep this
approximation controlled, the manifolds should belong to the same
algebraic family.

Generally speaking, such a sequence is hard to produce. For a K3 such
approximations are well known, and much used since the earliest works
on K3 in the 1960-es. It is known, in particular, that the variety of
quartic surfaces is dense in the moduli of all (non-algebraic) K3
surfaces.
 
In this paper, we generalize this theorem, proving that the moduli of
polarized hyperk\"ahler manifolds is a dense subset in the moduli of
all (non-algebraic) deformations. More precisely, given a rational
cohomology class  $\eta\in H^2(M,\mathbb Q)$, satisfying
$q(\eta,\eta)>0$, we show that any given $M$ can be approximated by
deformations of $M$ which satisfy $\eta\in H^{1,1}(M_1)$.

It is interesting that even this result seems to be quite hard to
prove. Our proof relies on rationality of $\eta$ and does not work when
$\eta$ is irrational,  though the statement is
most likely true in this case as well.

Another application of these results was obtained in 
\cite{_Kamenova_Verbitsky_}, a few years after the present
paper was finished. Jointly with Ljudmila Kamenova,
the second named author used the density theorem
to study the moduli space of hyperk\"ahler manifolds
admitting a Lagrangian fibration. It was shown that
this moduli space is a divisor, which is also dense in
the moduli of all deformations of the manifold. 
This theorem has some interesting further applications.
It was shown that the set of deformation classes of 
Lagrangian fibrations on a given hyperk\"ahler manifold
is finite. Also, it was shown that any hyperk\"ahler manifold
which has a deformation admitting a Lagrangian fibration
is Kobayashi non-hyperbolic. This proves, in particular,
that all known hyperk\"ahler manifolds are Kobayashi
non-hyperbolic, solving an old conjecture.

%%%%%%%%%%%%%%%%%%%%%%%%%%%%%%%%%%%%%%%%%%%%%%%%%%%%%%%%%%%%%%%
\subsection{Bogomolov-Beauville-Fujiki form and the mapping class
group}
\label{_BBF_Subsection_}
%%%%%%%%%%%%%%%%%%%%%%%%%%%%%%%%%%%%%%%%%%%%%%%%%%%%%%%%%%%%%%%

For a better understanding of the moduli space geometry, some basic
facts about topology of hyperk\"ahler manifolds should be stated. We
follow \cite{_V:Torelli_}.

\hfill

Let $\Omega$ be a holomorphic symplectic form on a hyperk\"ahler
manifold $M$. Bogomolov \cite{_Bogomolov:defo_} and Beauville
\cite{_Beauville_} defined the following bilinear symmetric $2$-form on
$H^2(M)$ :
\begin{equation}\label{_BBF_form_on_H^11_Equation_}
\begin{aligned}
\tilde q(\eta,\eta'):=&
\int_M\eta\wedge\eta'\wedge\Omega^{n-1}\wedge\bar\Omega^{n-1}-\\&
-\frac{(n-1)}n\frac{\Big(\int_M\eta\wedge\Omega^{n-1}\wedge\bar\Omega^n
\Big)\cdot\Big(\int_M \eta'\wedge\Omega^n\wedge\bar\Omega^{n-1}\Big)}
{\int_M\Omega^n\wedge\bar\Omega^n},
\end{aligned}
\end{equation}
where $4n=\dim_\mathbb RM$.

\hfill

The form $\tilde q$ is topological by its nature.

\hfill

%%%%%%%%%%%%%%%%%%%%%%%%%%%%%%%%%%%%%%%%%%%%%%%%
\alttheorem\label{_BBF_Theorem_}
\ \cite{_Fujiki:HK_}{\bf:}
{\sl Let\/ $M$ be a simple hyperk\"ahler manifold of real dimension\/
$4n$. Then there exist a bilinear, symmetric, primitive non-degenerate
integer\/ $2$-form\/
$q:H^2(M,\mathbb Z)\otimes H^2(M,\mathbb Z)\arrow\mathbb Z$ and a
positive constant\/ $c\in\mathbb Z$ such that\/
$\int_M\eta^{2n}=cq(\eta,\eta)^n$ for all\/ $\eta\in H^2(M)$. Moreover,
$q$ is proportional to the form\/ $\tilde q$ of\/
{\rm\eqref{_BBF_form_on_H^11_Equation_},} and has signature\/
$(3,b_2-3)$ {\rm(}with\/ $3$ pluses and\/ $b_2-3$ minuses\/{\rm).}}
\endproof

\hfill

Let $\Diff^+$ denote the group of orientation-preserving
diffeomorphisms of~$M$, and $\Diff^0$ its connected component, also
known as a group of isotopies. The quotient group
$\Gamma:=\Diff^+/\Diff^0$ is called the {\bf mapping class group} of
$M$. In \cite{_V:Torelli_} it was shown that $\Gamma$ preserves the
Bogomolov-Beauville-Fujiki form on $H^2(M)$ and that the corresponding
homomorphism to the orthogonal group
$\Gamma\arrow{\O}\big(H^2(M),q\big)$ has finite kernel. It was also
shown that the image of $\Gamma$ in ${\O}\big(H^2(M),q\big)$ is
commensurable to the group ${\O}\big(H^2(M,\mathbb Z),q\big)$ of
isometries of the integer lattice.

%%%%%%%%%%%%%%%%%%%%%%%%%%%%%%%%%%%%%%%%%%%%%%%%%%%%%%%%%%%%%%%
\subsection{Teichm\"uller space and the moduli space}
%%%%%%%%%%%%%%%%%%%%%%%%%%%%%%%%%%%%%%%%%%%%%%%%%%%%%%%%%%%%%%%

To state our main result in precise terms, we have to give a more
explicit description of the moduli space of a hyperk\"ahler manifold.
We follow \cite{_V:Torelli_}.

Let $M$ be a hyperk\"ahler manifold (compact and simple, as usual), and
$\Comp_0$ the Fr\`echet manifold of all complex structures of
hyperk\"ahler type on~$M$. The quotient $\Teich:=\Comp_0/\Diff^0$ of
$\Comp_0$ by isotopies is a finite-dimensional complex analytic space
by the same Kodaira-Spencer arguments as used to show that
$\Mod=\Comp/\Diff^+$ is complex analytic, where $\Comp$ is the Fr\`echet
manifold of all integrable complex, oriented structures on $M$. This
quotient is called {\bf the Teichm\"uller space} of $M$. When $M$ is a
complex curve, the quotient $\Comp/\Diff^0$ is the Teichm\"uller space
of this curve.

The mapping class group $\Gamma=\Diff^+/\Diff^0$ acts on $\Teich$ in the
usual way, and its quotient is the moduli space of $M$.

As shown in \cite{_Huybrechts:finiteness_}, $\Teich/\Gamma$ has a finite
number of connected components. Since $\Gamma$ is commensurable to 
a $SO(H^2(M,\Z), q)$, and $SO(H^2(M,\Z), q)$ acts virtually 
freely on $\Teich$ (\cite[Theorem 3.5, Theorem 4.29]{_V:Torelli_}),
the space $\Teich$ also has finitely many connected components.

Take a connected component $\Teich^I$
containing a given complex structure~$I$, and let
$\Gamma^I\subset\Gamma$ be the set of elements of $\Gamma$ fixing this
component. Since $\Teich$ has only a finite number of connected
components, $\Gamma^I$ has finite index in $\Gamma$. On the other hand,
as shown in \cite{_V:Torelli_}, the image of the group $\Gamma$ is
commensurable to ${\O}\big(H^2(M,\mathbb Z),q\big)$.

In \cite[Lemma 2.6]{_V:Torelli_} it was proved that any hyperk\"ahler
structure on a given simple hyperk\"ahler manifold is also simple.
Therefore, $H^{2,0}(M,I')=\mathbb C$ for all $I'\in\Comp$. This
observation is a key to the following well-known definition.

\hfill

%%%%%%%%%%%%%%%%%%%%%%%%%%%%%%%%%%%%%%%%%%%%%%%%
\definition
Let $(M,I)$ be a hyperk\"ahler manifold, and $\Teich$ its Teichm\"uller
space. Consider a map $\Per:\Teich\arrow\mathbb PH^2(M,\mathbb C)$,
sending $J$ to the line $H^{2,0}(M,J)\in\mathbb PH^2(M,\mathbb C)$. It
is easy to see that $\Per$ maps $\Teich$ into the open subset of a
quadric, defined by
\begin{equation*}\label{_perspace_lines_Equation_}
\Perspace:=\big\{l\in\mathbb PH^2(M,\mathbb C)\ \big|\ q(l,l)=0,\
q(l,\bar l)>0\big\}.
\end{equation*}
The map $\Per:\Teich\arrow\Perspace$ is called the {\bf period map},
and the set $\Perspace$ the {\bf period space}.

\hfill

The following fundamental theorem is due to F. Bogomolov
\cite{_Bogomolov:defo_}.

\hfill

%%%%%%%%%%%%%%%%%%%%%%%%%%%%%%%%%%%%%%%%%%%%%%%%%%%%%%%%%%%%
\theorem
{\sl Let\/ $M$ be a simple hyperk\"ahler manifold, and $\Teich$ its
Teichm\"uller space. Then the period map\/ $\Per:\Teich\arrow\Perspace$
is a local diffeomorphism\/ {\rm(}that is, an etale map\/{\rm).}
Moreover, it is holomorphic.}
\endproof

\hfill

%%%%%%%%%%%%%%%%%%%%%%%%%%%%%%%%%%%%%%%%%%%%%%%%
\remark
{\sl Bogomolov's theorem implies that $\Teich$ is smooth. However, it
is not necessarily Hausdorff\/ {\rm(}and it is non-Hausdorff even in
the simplest examples\/{\rm).}}
\endproof

%%%%%%%%%%%%%%%%%%%%%%%%%%%%%%%%%%%%%%%%%%%%%%%%%%%%%%%%%%%%%%%
\subsection{The polarized Teichm\"uller space}
%%%%%%%%%%%%%%%%%%%%%%%%%%%%%%%%%%%%%%%%%%%%%%%%%%%%%%%%%%%%%%%

In \cite[Corollary 2.6]{_Verbitsky:parabolic_}, the following
proposition was deduced from \cite{_Boucksom_}
and~\cite{_Demailly_Paun_}.

\hfill

%%%%%%%%%%%%%%%%%%%%%%%%%%%%%%%%%%%%%%%%%%%
\theorem\label{_Kahler_cone_Pic=1_Theorem_}
{\sl Let\/ $M$ be a simple hyperk\"ahler manifold, such that all
integer\/ $(1,1)$-classes satisfy\/ $q(\nu,\nu)\geq 0$. Then its
K\"ahler cone is one of two connected components of the set\/
$K:=\big\{\nu\in H^{1,1}(M,\mathbb R)\ \big|\ q(\nu,\nu)>0\big\}$.}
\endproof

\hfill

Consider an integer vector $\eta\in H^2(M)$ which is positive, that is,
satisfies $q(\eta, \eta)>0$. Denote by $\Teich^\eta$ the set of all
$I\in \Teich$ such that $\eta$ is of type $(1,1)$ on $(M, I)$. The
space $\Teich^\eta$ is a closed divisor in $\Teich$. Indeed, by
Bogomolov's theorem, the period map $\Per:\Teich\arrow\Perspace$ is
etale, but~the image of $\Teich^\eta$ is the set of all $l\in\Perspace$
which are orthogonal to $\eta$; this condition defines a closed divisor
$C_\eta$ in $\Perspace$, hence $\Teich^\eta=\Per^{-1}(C_\eta)$ is also
a closed divisor.

When $I\in\Teich^\eta$ is generic, Bogomolov's theorem implies that the
space of rational $(1,1)$-classes $H^{1,1}(M,\mathbb Q)$ is
one-dimensional and generated by $\eta$. This is seen from the
following argument.  Locally around a given point $I$ the period map
$\Teich^\eta\arrow\Perspace$ is surjective on the set $\Perspace^\eta$
of all $I\in\Perspace$ for which $\eta\in H^{1,1}(M,I)$. However, the
Hodge-Riemann relations give 
\begin{equation}\label{_Perspace^eta_ortho_Equation_}
\Perspace^\eta=\big\{l\in\Perspace\ \big|\ q(\eta,l)=0\big\}.
\end{equation}
Denote the set of such points of $\Teich^\eta$ by
$\Teich^\eta_{\text{gen}}$. It follows from
\ref{_Kahler_cone_Pic=1_Theorem_} that, for any
$I\in\Teich^\eta_{\text{gen}}$, either $\eta$ or $-\eta$ is a K\"ahler
class on $(M,I)$.

Consider a connected component $\Teich^{\eta,I}$ of $\Teich^\eta$.
Changing the sign of $\eta$ if necessary, we may assume that $\eta$ is
K\"ahler on $(M,I)$. By Kodaira's theorem about stability of K\"ahler
classes, $\eta$ is K\"ahler in some neighbourhood
$U\subset\Teich^{\eta,I}$ of $I$. Therefore, the sets
$$V_+:=\big\{I\in\Teich^\eta_{\text{gen}}\ \big|\ \eta\text{ is
K\"ahler on }(M,I)\big\}$$
and
$$V_-:=\big\{I\in\Teich^\eta_{\text{gen}}\ \big|\ -\eta\text{ is
K\"ahler on }(M,I)\big\}$$ 
are open in $\Teich^\eta_{\text{gen}}$. It is easy to see that
$\Teich^\eta_{\text{gen}}$ is a complement to a union of countably many
divisors in $\Teich^\eta$ corresponding to the points
$I'\in\Teich^\eta$ with $\rk\Pic(M,I')>1$. Therefore, for any connected
open subset $U\subset\Teich^\eta$, the intersection
$U\cap\Teich^\eta_{\text{gen}}$ is connected. Since
$\Teich^\eta_{\text{gen}}$ is represented as a disjoint union of open
sets $V_+\sqcup V_-$, every connected component 
of $\Teich^\eta$ is contained in $V_+$
or in $V_-$. We obtained the following corollary.

\hfill

%%%%%%%%%%%%%%%%%%%%%%%%%%%%%%%%%%%%%%
\corollary
{\sl Let\/ $\eta\in H^2(M)$ be a positive integer vector, $\Teich^\eta$
the corresponding divisor in the Teichm\"uller space, and\/
$\Teich^{\eta,I}$ a connected component of\/ $\Teich^\eta$ containing a
complex structure\/ $I$. Assume that\/ $\eta$ is K\"ahler on\/ $(M,I)$.
Then\/ $\eta$ is K\"ahler for all\/ $I'\in\Teich^{\eta,I}$ which
satisfy\/ $\rk H^{1,1}(M,\mathbb Q)=1$.}
\endproof

\hfill

We call the set $\Teich^\eta_{\text{pol}}$ of all $I\in\Teich^\eta$ for
which $\eta$ is K\"ahler the {\bf polarized Teichm\"uller space}, and
$\eta$ its {\bf polarization}. From the above arguments it is clear
that the polarized Teichm\"uller space is open and dense in
$\Teich^\eta$.

The quotient ${\cal M}_\eta$ of $\Teich^\eta_{\text{pol}}$ by the
subgroup of a mapping class group fixing $\eta$ is called the {\bf
moduli of polarized hyperk\"ahler manifolds}. It is known (due to the
general theory which goes back to Viehweg and Grothen\-dieck that
${\cal M}_\eta$ is Hausdorff and quasiprojective (see
e.g.~\cite{_Viehweg:moduli_} and \cite{_GHK:moduli_HK_}).

We conclude that there are countably many quasiprojective divisors
${\cal M}_\eta$ immersed in the moduli space $\Mod$ of hyperk\"ahler
manifolds. Moreover, every algebraic complex structure belongs to one
of these divisors. However, these divisors need not to be closed.
Indeed, as we prove in this paper, each of the ${\cal M}_\eta$ is dense
in $\Mod$.

The main result of the present paper is the following theorem.

\hfill

%%%%%%%%%%%%%%%%%%%%%%%%%%%%%%%%%%%%%%%%%%%%%%%%%%%%%%%%%%%%%
\theorem\label{_dense_main_Theorem_}
{\sl Let\/ $M$ be a compact, simple hyperk\"ahler manifold, $\Teich^I$
a connected component of its Teichm\"uller space, and\/
$\Teich^I\stackrel\Psi\arrow\Teich^I/\Gamma^I=\Mod$ its projection to
the moduli of complex structures. Consider a positive vector\/
$\eta\in H^2(M,\mathbb Z)$, and let\/ $\Teich^{I,\eta}$ be the
corresponding connected component of the polarized Teichm\"uller space.
Assume that\/ $b_2(M)>3$. Then the image\/ $\Psi(\Teich^{I,\eta})$ is
dense in\/ $\Mod$.}

\hfill

We deduce \ref{_dense_main_Theorem_} from
\ref{_lattices_dense_Proposition_} in Section \ref{_Torelli_Section_},
and prove \ref{_lattices_dense_Proposition_} in Section
\ref{_lattices_Section_}.

\hfill

\remark
{\sl We assumed positivity of $\eta$ in the statement of
\ref{_dense_main_Theorem_}, but this assumption is completely
unnecessary. In fact, for\/ $\eta$ non-positive, the proof of
\ref{_dense_main_Theorem_} becomes easier\/
{\rm(}\ref{negative_Remark}{\rm).}}
\endproof

%%%%%%%%%%%%%%%%%%%%%%%%%%%%%%%%%%%%%%%%%%%%%%%%%%%%%%%%%%%%%%%
\section{Torelli theorem and polarizations}
\label{_Torelli_Section_}
%%%%%%%%%%%%%%%%%%%%%%%%%%%%%%%%%%%%%%%%%%%%%%%%%%%%%%%%%%%%%%%

In this Section, we reduce \ref{_dense_main_Theorem_} to a statement
about lattices and arithmetic groups, proven in Section
\ref{_lattices_Section_}.

Let $M$ be a topological space, not necessarily Hausdorff. We say that
points $x,y\in M$ are {\bf inseparable} (denoted $x\sim y$) if for any
open subsets $U\ni x,V\ni y$, one has $U\cap V\ne\varnothing$. 

\hfill

%%%%%%%%%%%%%%%%%%%%%%%%%%%%%%%%%%%%%%%%%%%%%%%%
\alttheorem\label{_Torelli_Theorem_}
\cite[Theorem 1.14, Theorem 1.16]{_V:Torelli_}{\bf:}
{\sl Let\/ $\Teich$ be a Teichm\"uller space of a hyperk\"ahler
manifold, and\/ $\sim$ the inseparability relation defined above.
Then\/ $\sim$ is an equivalence relation, and the quotient\/
$\Teich_b:={\Teich}/{\sim}$ is a smooth, Hausdorff, complex analytic
manifold. Moreover, the period map\/ $\Per:\Teich\arrow\Perspace$
induces a complex analytic diffeomorphism\/ $\Teich_b^I\arrow\Perspace$
for each connected component\/ $\Teich_b^I$ of\/ $\Teich_b$.}~\endproof

\hfill

%%%%%%%%%%%%%%%%%%%%%%%%%%%%%%%%%%%%%%%%%%%%%%%%
\remark 
{\sl As shown by Huybrechts\/ {\rm\cite{_Huybrechts:basic_},}
inseparable points on a Teichm\"uller space correspond to
bimeromorphically equivalent hyperk\"ahler\break manifolds. The
Hausdorff quotient\/ $\Teich_b={\Teich}/{\sim}$ is called the {\bf
birational Teichm\"uller space} of\/ $M$.}
\endproof

\hfill

By construction, the action of the mapping class group $\Gamma$ on
$\Teich_b$ is compatible with the natural action of
${\O}\big( H^2(M,\mathbb Z),q\big)$ on $\Perspace$. Define the {\bf
birational moduli space} as $\Mod_b:=\Teich_b/\Gamma$. The space
$\Mod_b$ is obtained by gluing together some (not all) inseparable
points in $\Mod$. By~\ref{_Torelli_Theorem_},
$\Mod_b=\Perspace/\Gamma^I$, where $\Gamma^I$ is a subgroup of $\Gamma$
fixing a connected component $\Teich^I$ of the Teichm\"uller space. As
follows from \cite[Theorem 3.5]{_V:Torelli_} (see also Subsection
\ref{_BBF_Subsection_}), the image of $\Gamma^I$ in $\Aut(\Perspace)$
is a finite index subgroup in ${\O}\big(H^2(M,\mathbb Z),q\big)$.

It is well known that the homogeneous space 
$$\Perspace=\big\{l\in\mathbb PH^2(M,\mathbb C)\ \big|\ q(l,l)=0,\
q(l,\bar l)>0\big\}$$
is naturally identified with the Grassmanian
$${\Gr}^{++}\big(H^2(M,\mathbb
R)\big)\cong\SO(3,b_2-3)/\SO(2)\times\SO(1,b_2-3)$$
of oriented positive 2-dimensional planes in $H^2(M,\mathbb R)$. This
identification is performed as follows: to each line
$l\in\mathbb PH^2(M,\mathbb C)$ one associates the plane spanned by
$\Re(l),\Im(l)$. Under this identification, the image of the polarized
Teichm\"uller space $\Teich^\eta$ is the space of all $2$-dimensional
planes $P\in\Gr^{++}\big(H^2(M,\mathbb R)\big)$ orthogonal to $\eta$
(see \eqref{_Perspace^eta_ortho_Equation_}). Then
\ref{_dense_main_Theorem_} is implied by the following statement.

\hfill

%%%%%%%%%%%%%%%%%%%%%%%%%%%%%%%%%%%%%%%%%%%%%%%%%%%
\theorem\label{_Densiti_in_Gr_cohomo_Theorem_}
{\sl Let\/ $M$ be a simple, compact hyperk\"ahler manifold,
$V:=H^2(M,\mathbb R)$ its second cohomology, $L:=H^2(M,\mathbb Z)$,
and\/ $q$ the Bogomolov-Beauville-Fujiki form on\/ $V$. Given a
positive integer vector\/ $\eta\in L$, denote by\/
$\Gr^{++}(\eta^\bot)\subset\Gr^{++}(V)$ the space of all planes
orthogonal to\/ $\eta$. Consider a finite index subgroup\/
$G\subset\SO\big(H^2(M,\mathbb Z),q\big)$ acting on\/ $\Gr^{++}(V)$ in
the natural way. Then\/ $G\cdot\Gr^{++}(\eta^\bot)$ is dense in\/
$\Gr^{++}(V)=\Perspace$.}

\hfill

\ref{_Densiti_in_Gr_cohomo_Theorem_} is implied by a more general
\ref{_lattices_dense_Proposition_} proven in the next section using the
framework laid down in \cite{_AG:Geometries_}.

\hfill

%%%%%%%%%%%%%%%%%%%%%%%%%%%%%%%%%%%%%%%%%%%%%%%%
\remark \label{_K3_Remark_}
{\sl When\/ $M$ is a\/ {\rm K3} surface, the Bogomolov-Beauville-Fujiki
form is unimodular, and the mapping class group is generated by
appropriate reflections. From a statement known as ``Eichler's
criterion''\/ {\rm(}see\/
{\rm\cite[Proposition 3.3(i)]{_GHK:modular_}{\rm),}} the mapping class
group acts transitively on the set of integer vectors of a given length
in\/ $H^2(M)$. \ref{_Densiti_in_Gr_cohomo_Theorem_} follows from this
observation easily. When the Eichler's criterion cannot be applied, its
proof is more complicated.}
\endproof

%%%%%%%%%%%%%%%%%%%%%%%%%%%%%%%%%%%%%%%%%%%%%%%%
\section{Arithmetic subgroups in $\O(p,q)$}
\label{_lattices_Section_}
%%%%%%%%%%%%%%%%%%%%%%%%%%%%%%%%%%%%%%%%%%%%%%%%

Let $V$ be a finite-dimensional $\mathbb R$-vector space equipped with
a non-degenerate symmetric form $\langle\cdot,\cdot\rangle$ and $W$ an
$\mathbb R$-vector subspace in $V$. Denote by $\Gr_{++}(W)$
(respectively, by $\Gr_{+-}(W)$) the part of the Grassmannian
$\Gr_\mathbb R(2,V)$ of $2$-dimensional $\mathbb R$-subspaces in $V$
formed by the subspaces of signature $++$ (respectively, $+-$) in $W$.

\hfill

%%%%%%%%%%%%%%%%%%%%%%%%%%%%%%%%%%%%%%%%%%%%%%%%
\definition
We shall call a discrete, additive subgroup $L\subset V$ a {\bf
lattice} if $V=\mathbb R\otimes_\mathbb ZL$ and
$\langle l_1,l_2\rangle\in\mathbb Q$ for all $l_1,l_2\in L$. Denote by
$\O(V)$ and $\O(L)$ the corresponding orthogonal groups:
\begin{align*}
\O(V):=&\Big\{g\in\GL(V)\ \Big|\ \big\langle
g(v_1),g(v_2)\big\rangle=\langle v_1,v_2\rangle\text{ for all
}v_1,v_2\in V\Big\},\\
\O(L):=&\big\{g\in\O(V)\ \big|\ g(L)=L\big\}.
\end{align*}
Clearly, $\O(V)$ acts on $\Gr_{++}(V)$. For $S\subset V$, we denote
$$S^\perp:=\big\{v\in V\ \big|\ \langle v,S\rangle=0\big\}.$$

The purpose of the present section is to prove

\hfill

%%%%%%%%%%%%%%%%%%%%%%%%%%%%%%%%%%%%%%%%%%%%%%%%
\proposition\label{_lattices_dense_Proposition_}
{\sl Let\/ $V$ be an\/ $\mathbb R$-vector space equipped with a
non-degene\-rate symmetric form of signature\/ $(s_+,s_-)$ with\/
$s_+\ge3$ and\/ $s_-\ge1$. Consider a lattice\/ $L\subset V$. Let\/
$\Gamma$ be a subgroup of finite index in\/ $\O(L)$, and\/ $l\in L$ a
{\bf positive vector}, i.e., one which satisfies\/
$\langle l,l\rangle>0$. Then\/ $\Gamma\cdot\Gr_{++}(l^\perp)$ is dense
in $\Gr_{++}(V)$.}

\hfill

The proof of \ref{_lattices_dense_Proposition_} takes the rest of this
Section.

\hfill

{\bf Proof of \ref{_lattices_dense_Proposition_}: Step 1:}
We reduce \ref{_lattices_dense_Proposition_} to a case of a space $V$
of signature $(3,1)$.

A subspace $W\subset V$ is called {\bf rational} if
$\rk(W\cap L)=\dim_\mathbb RW$ or, equivalently, if $W=\mathbb RW_0$
with a $\mathbb Q$-subspace $W_0\subset\mathbb QL$. Since the rational
subspaces are dense in $\Gr_{++}(V)$, it suffices to show that an
arbitrary rational $2$-plane $C\in\Gr_{++}(V)$ belongs to the closure
of $\Gamma\cdot\Gr_{++}(l^\perp)$. We have $C=\mathbb RC_0$ for some
$\mathbb Q$-subspace $C_0\subset\mathbb QL$.

Obviously, $\mathbb QL$ has signature $(s_+,s_-)$. Applying to
$\mathbb Q$-subspaces in $\mathbb QL$ the standard orthogonalization
arguments, we can find a $\mathbb Q$-subspace $U_0\subset\mathbb QL$ of
signature $+++-$ that contains both $l$ and $C_0$. Indeed, we have
$l=c_0+c_1$, where $c_0\in C_0$ and $c_1\in\mathbb QL\cap C_0^\perp$
with $\mathbb QL\cap C_0^\perp$ of signature $(s_+-2,s_-)$. We can
always pick a $2$-dimensional $\mathbb Q$-subspace
$C_1\subset\mathbb QL\cap C_0^\perp$ of signature $+-$ that includes
$c_1$ and put $U_0:=C_0\oplus C_1$. So, $U:=\mathbb RU_0$ is rational
of signature $+++-$ and $L_0:=U\cap L$ is a lattice in $U$. To prove
that the $2$-plane $C$ belongs to the closure of
$\Gamma\cdot\Gr_{++}(l^\perp)$, it would suffice to show that the set
$(\Gamma\cap\Gamma')\cdot\Gr_{++}(l^\perp\cap U)$ is dense in the
corresponding $++$-Grassmannian $\Gr_{++}(U)$, where
$\Gamma':=\O(L_0)$.

\hfill

{\bf Step 2:}
We prove that the orthogonal groups $\O(L)$ and $\O(L')$ are {\bf
commensurable}, i.e., the subgroup $\O(L)\cap\O(L')$ has finite index
in $\O(L)$ and in $\O(L')$, if lattices $L,L'\subset V$ are
commensurable.

Taking $L\cap L'$ for $L'$, we can assume that $mL\subset L'\subset L$
for some $0\ne m\in\mathbb Z$. Put
$\overline{L'}:=L'/mL\subset\overline L:=L/mL$ and note that $\O(L)$
acts on $\overline L$ because $\O(L)=\O(mL)$. We can see that the group
$\O(L)\cap\O(L')=\big\{g\in\O(L)\ \big|\ g(L')=L'\big\}$ coincides with
the stabilizer $\St_{\O(L)}\overline{L'}$. Hence, $\O(L)\cap\O(L')$ has
finite index in $\O(L)$. Since $m(\frac1mL')\subset L\subset\frac1mL'$
and $\O(\frac1mL')=\O(L')$, we infer as well that $\O(L)\cap\O(L')$ has
finite index in~$\O(L')$.

\hfill

{\bf Step 3:}
Let $W\subset V$ be a rational non-degenerate subspace. Then we have an
orthogonal decomposition $V=W\oplus W^\perp$ and $W^\perp$ is rational.
Define $L_0:=W\cap L$, $L_1:=W^\perp\cap L$, and $L':=L_0+L_1$. It is
immediate that $L'$ is a lattice in $V$ such that
$\mathbb QL'=\mathbb QL$. By Step 2, the orthogonal groups $\O(L)$ and
$\O(L')$ are commensurable. Since $\O(L_0)\times\O(L_1)\subset\O(L')$,
there exists a subgroup $\Gamma_0$ of finite index in $\O(L_0)$ such
that $\Gamma_0\subset\Gamma$.

\hfill

{\bf Step 4:}
We reduce \ref{_lattices_dense_Proposition_} to
\ref{_3-dim-reformula_Lemma_} below.

Applying Steps 1 and 3, we can assume that $(s_+,s_-)=(3,1)$. Indeed,
by Step 1, we need only to show that
$(\Gamma\cap\Gamma')\cdot\Gr_{++}(l^\perp\cap U)$ is dense in
$\Gr_{++}(U)$, where $\Gamma':=\O(L_0)$, $L_0:=U\cap L$, and
$U\subset V$ is a rational subspace of signature $+++-$. Taking $W:=U$
in Step 3, we find a subgroup $\Gamma_0$ of finite index in $\Gamma'$
such that $\Gamma_0\subset\Gamma$.

Now using the homeomorphism $\Gr_{++}(V)\to\Gr_{+-}(V)$,
$G\mapsto G^\perp$, i.e., taking instead of subspaces of signature
$++$, their orthogonal complements (of signature $+-$), we reformulate
\ref{_lattices_dense_Proposition_} as follows:

\hfill

\centerline{\sl Every rational\/ $G_0\in\Gr_{+-}(V)$ belongs to}

\centerline{\sl the closure of\/
$\Gamma\cdot\big\{G\in\Gr_{+-}(V)\ \big|\ G\ni l\big\}$.}

\hfill

The subspace $W$ spanned by $l,G_0$ is rational of signature $++-$.
Again using Step 3, we reduce \ref{_lattices_dense_Proposition_} to

\hfill

%%%%%%%%%%%%%%%%%%%%%%%%%%%%%%%%%%%%%%%%%%%%%%%%%%%
\lemma\label{_3-dim-reformula_Lemma_}
{\sl Let\/ $V$ be an\/ $\mathbb R$-vector space equipped with a
symmetric form of signature\/ $++-$, $\Gamma$ a subgroup of finite
index in\/ $\O(L)$, where\/ $L$ is a lattice in\/ $V$, and\/ $l\in V$ a
positive vector. Then\/
$\Gamma\cdot\big\{G\in\Gr_{+-}(V)\ \big|\ G\ni l\big\}$ is dense in\/
$\Gr_{+-}(V)$.}

\hfill

Till the end of this Section we fix $\Gamma$ as in
\ref{_3-dim-reformula_Lemma_}.

\hfill

In fact, we deal now with a hyperbolic plane
$\overline{\mathbb H}_\mathbb R^2=\mathbb H_\mathbb
R^2\sqcup\partial\mathbb H_\mathbb R^2$.
Let us state in \ref{geodesics_}, \ref{orthogonal_geodesics}, and
\ref{isometries_} a few simple and well-known facts concerning the
hyperbolic plane (see e.g.~\cite{_AG:Geometries_}).

\hfill

\claim\label{geodesics_}
{\sl The plane\/ $\overline{\mathbb H}_\mathbb R^2$ can be identified
with the set of all nonpositive points in the real projective plane\/
$\mathbb P_\mathbb RV$, where the isotropic ones form the {\bf
absolute} $\partial\mathbb H_\mathbb R^2$. In the affine chart related
to orthonormal coordinates on\/~$V$, the plane\/
$\overline{\mathbb H}_\mathbb R^2$ is nothing but a closed unitary
disc. In this way, we~obtain the Beltrami-Klein model of a hyperbolic
plane, where geodesics are chords of the disc. In other words, we can
describe a geodesic in\/ $\overline{\mathbb H}_\mathbb R^2$ as the
projectivization\/
$\mathbb P_\mathbb RG\cap\overline{\mathbb H}_\mathbb R^2$ of a
subspace\/ $G\in\Gr_{+-}(V)$. We keep denoting this geodesic by\/ $G$.
Of course, every geodesic\/ $G$ can be described via its {\bf vertices}
$v,v'\in\partial\mathbb H_\mathbb R^2$ as\/ $G=[v,v']$. In terms of\/
$V$, this means that the\/ $\mathbb R$-vector subspace\/ $G$ is spanned
by\/ $v,v'$.}
\endproof

\hfill

\claim\label{orthogonal_geodesics}
{\sl Let\/ $G'\subset\overline{\mathbb H}_\mathbb R^2$ be a geodesic
not passing through a point\/ $v\in\partial\mathbb H_\mathbb R^2$,
i.e., $v\notin G'$. Then, reflecting\/ $v$ in\/ $G'$, we obtain a
point\/ $v'\in\partial\mathbb H_\mathbb R^2$ such that the geodesics\/
$G'$ and\/ $[v,v']$ are orthogonal.}
\endproof

\hfill

\claim\label{isometries_}
{\sl The group\/ $\O(V)$ acts naturally on\/
$\overline{\mathbb H}_\mathbb R^2$. On\/ $\mathbb H_\mathbb R^2$, the
group\/ $\O(V)$ acts by isometries.}
\endproof

\hfill

We can now reduce \ref{_3-dim-reformula_Lemma_} to the following
statement about the hyperbolic plane:

\hfill

%%%%%%%%%%%%%%%%%%%%%%%%%%%%%%%%%%%%%%%%%%%%%%%%%%%%%%%
\lemma\label{_geode_dense_Lemma_}
{\sl Let\/ $G'$ be a geodesic on the hyperbolic plane\/
$\mathbb H_\mathbb R^2$, and\/ $\Gamma\cdot G'$ the set of all
geodesics obtained from\/ $G'$ by the action of\/ $\Gamma$. Then the
set of all geodesics orthogonal to some\/ $G''\in\Gamma\cdot G'$ is
dense in the set of all geodesics in\/ $\mathbb H_\mathbb R^2$.}

\hfill

{\bf Reduction of \ref{_3-dim-reformula_Lemma_} to 
\ref{_geode_dense_Lemma_}.}
Let $G'$ be the orthogonal complement of
$l\in\mathbb H_\mathbb R^2\subset {\mathbb P}V$, considered as a
geodesic in $\mathbb H_\mathbb R^2$. It is easy to see that the
inclusion $G\ni l$ is equivalent to the fact that the geodesics $G$ and
$G':=l^\perp$ are orthogonal (see, for instance, the duality described
in the introductory \cite[Section 1]{_AG:Geometries_} shortly after
Example 1.7). For this choice of~$G'$, \ref{_3-dim-reformula_Lemma_} is
clearly equivalent to \ref{_geode_dense_Lemma_}.
\endproof

\hfill

We reduce \ref{_geode_dense_Lemma_} further, obtaining a simpler
statement about the hyperbolic plane:

\hfill

\lemma\label{hyperbolic_Lemma}
{\sl Let\/ $v,v'\in\partial\mathbb H_\mathbb R^2$ be distinct points on
the absolute and\/ $G'$ a geodesic. For every $\gamma\in\Gamma$, denote
by\/ $R_\gamma$ the reflection in the geodesic\/ $\gamma(G')$. Then\/
$v'$ belongs to the closure of the set\/
$\big\{R_\gamma(v)\ \big|\ \gamma\in\Gamma,\ v\notin\gamma(G')\big\}$
formed by the reflections of\/ $v$ in those geodesics\/ $\gamma(G')$
that do not pass through\/ $v$.}

\hfill

{\bf Reduction of \ref{_geode_dense_Lemma_} to \ref{hyperbolic_Lemma}:}
By \ref{orthogonal_geodesics}, the geodesic $[v,R_\gamma(v)]$ is
orthogonal to $\gamma(G')$. To prove \ref{_geode_dense_Lemma_}, it
suffices to show that the set of such geodesics is dense in the set of
all geodesics of the form $[v,v']$, where $v$ is fixed.
\ref{hyperbolic_Lemma} says that we are able to approximate $v'$ by
$R_\gamma(v)$ for an appropriate $\gamma\in\Gamma$. Hence, we can
approximate the geodesic $[v,v']$ by geodesics orthogonal to some
$G''\in\Gamma\cdot G'$.
\endproof

\hfill

We deduce \ref{hyperbolic_Lemma} from two easy lemmas below,
\ref{fixed_points_density_Lemma} and \ref{no_fix_point_Lemma}. First,
we need a few more simple and well-known facts concerning the
hyperbolic plane:

\hfill

\begin{wrapfigure}[10]{r}{0.22\textwidth}
\includegraphics[width=0.22\textwidth]{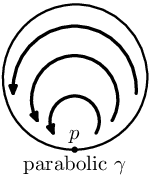}
\end{wrapfigure}
$\bullet$ The nontrivial orientation-preserving isometries of
$\mathbb H_\mathbb R^2$ are classified with respect to the location of
their fixed points: an {\bf elliptic} one has a (unique) fixed point in
$\mathbb H_\mathbb R^2$; a {\bf hyperbolic} one has exactly two fixed
points on the absolute; and a {\bf parabolic} one has exactly one fixed
point on the absolute.

\hfill

$\bullet$ Let $p\in\partial\mathbb H_\mathbb R^2$ be the fixed point of
a parabolic isometry $\gamma$ and let
$v\in\partial\mathbb H_\mathbb R^2$. Then $\gamma^n(v)\to p$ as
$n\to\infty$.

\hfill

\begin{wrapfigure}[10]{r}{0.27\textwidth}
\includegraphics[width=0.27\textwidth]{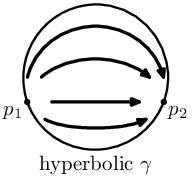}
\end{wrapfigure}
$\bullet$ The fixed points in $\partial\mathbb H_\mathbb R^2$ of a
hyperbolic isometry $\gamma$ are the {\bf repeller} $p_1$ and the {\bf
attractor} $p_2$. This means that, for every
$v\in\partial\mathbb H_\mathbb R^2$ such that $v\ne p_1$, we have
$\gamma^n(v)\to p_2$ as $n\to\infty$. When taking $\gamma^{-1}$ in
place of $\gamma$, the repeller becomes the attractor and vice versa.

\hfill

We arrive at the following remark needed in the proof of
\ref{hyperbolic_Lemma}.

\hfill

\remark\label{limit_Remark}
{\sl Let $\gamma$ be a hyperbolic or parabolic isometry,
$p\in\partial\mathbb H_\mathbb R^2$ a fixed point of $\gamma$, and
$u,u'\in\partial\mathbb H_\mathbb R^2$ points not fixed by $\gamma$.
Then, for $n\to\infty$ or for $n\to-\infty$, both limits
$\lim\gamma^n(u)$ and $\lim\gamma^n(u')$ exist and are equal to $p$.}
\endproof

\hfill

\lemma\label{fixed_points_density_Lemma}
{\sl The set\/
$F:=\big\{p\in\partial\mathbb H_\mathbb R^2\ \big|\ \gamma(p)=p\text{
\rm for some }1\ne\gamma\in\Gamma\big\}$
of points on the absolute fixed by some nontrivial\/ $\gamma\in\Gamma$
is dense in\/ $\partial\mathbb H_\mathbb R^2$.}

\hfill

{\bf Proof:}
Suppose that there exists an open arc
$A\subset\partial\mathbb H_\mathbb R^2$ such that
$A\cap F=\varnothing$. By the Zorn lemma, we can take maximal $A$ with
this property. By~construction, $\gamma(A)$ also enjoys the property of
the maximality for every $\gamma\in\Gamma$. Every point on the boundary
$\partial A$ belongs to the closure of $F$. Let
$\gamma\in\Gamma$. Then $A\cap\gamma(A)=\varnothing$ or $A=\gamma(A)$
because otherwise $\gamma(A)$ contains an open neighbourhood of one end
of $A$, which intersects $F$.

Due to B.~A.~Venkov (see \cite[Example 7.5, p.~33]{_VGSh:VINITI_}),
$\O(L)$ is known to act discretely on $\mathbb H_\mathbb R^2$, is
finitely generated, and is of finite coarea. Note that Selberg's
Theorem \cite[Theorem 3.2, p.~18]{_VGSh:VINITI_} claims that every
finitely generated matrix group over a field of characteristic $0$ has
a subgroup of finite index without torsion. Therefore, we can at the
very beginning pass to a torsion-free subgroup of finite index in
$\Gamma$ thus assuming that all isometries in $\Gamma$ are
orientation-preserving and that there are no elliptic isometries in
$\Gamma$.

Let $\partial A=\{p,p'\}$. Since $\Gamma$ has no elliptic isometries
and all isometries in $\Gamma$ are orientation-preserving, the
stabilizer $\Gamma'':=\St_\Gamma A$ of $A$ in $\Gamma$ is a discrete
\begin{wrapfigure}[11]{l}{0.36\textwidth}
\includegraphics[clip=true,trim=0 0 0
2.1,width=0.36\textwidth]{pictures-6.eps}
\end{wrapfigure}
group of orientation-preserving isometries of the geodesic $[p,p']$.
Hence, $\Gamma''$ is cyclic, generated by some $\gamma_0\ne1$. Let $G$
be a geodesic perpendicular to $[p,p']$. Consider the open region
$D\subset\mathbb H_\mathbb R^2$ limited by
$A\cup G\cup[p,p']\cup\gamma_0(G)$. It is easy to see that
$D\cap\gamma(D)=\varnothing$ for every $1\ne\gamma\in\Gamma''$. For
any $\gamma\in\Gamma\setminus\Gamma''$, we have
$A\cap\gamma(A)=\varnothing$, which again implies
$D\cap\gamma(D)=\varnothing$. Therefore, $D$ is a part of a fundamental
domain for $\Gamma$. Since the area of $D$ is infinite, we arrive at a
contradiction.
\endproof

\hfill

\lemma\label{no_fix_point_Lemma}
{\sl Let\/ $u,u'\in\partial\mathbb H_\mathbb R^2$ be distinct points.
Then there exists a hyperbolic or parabolic\/ $\gamma_0\in\Gamma$ such
that\/ $\gamma_0(u)\ne u$ and\/ $\gamma_0(u')\ne u'$.}

\hfill

{\bf Proof:}
As in the proof of \ref{fixed_points_density_Lemma}, we assume $\Gamma$
torsion-free. Suppose that $\gamma_0(u)=u$ or $\gamma_0(u')=u'$ for
every $\gamma_0\in\Gamma$. If $\gamma,\gamma'\in\Gamma$ fix
respectively $u,u'$ and do not fix respectively $u',u$, then
$\gamma\gamma'$ does not fix both $u$ and $u'$. Therefore, we can
assume that $\gamma(u)=u$ for all $\gamma\in\Gamma$. It is well known
(consider the upper half-plane model with $u=0$) that the group of all
orientation-preserving isometries of $\mathbb H_\mathbb R^2$ is
isomorphic to $\PSL_2(\mathbb R)$ and that
$S:=\St_{\PSL_2(\mathbb R)}u\simeq\left\{\left[\smallmatrix\alpha&0\\
a&\alpha^{-1}\endsmallmatrix\right]\Big|\ a,\alpha\in\mathbb R,\
\alpha>0\right\}$.
Since $S$ is Zariski closed, the inclusion $\Gamma\subset S$ would
contradict the Borel density theorem
\cite[Theorem 8.2, p.~37]{_VGSh:VINITI_} which implies that $\Gamma$
should be Zariski dense in $\PSL_2(\mathbb R)$.
\endproof

\hfill

{\bf Proof of \ref{hyperbolic_Lemma}:}
For suitable distinct points $u,u'\in\partial\mathbb H_\mathbb R^2$,
the geodesic $G'$ in \ref{hyperbolic_Lemma} has the form $G'=[u,u']$.

Let $A$ be a small connected open neighbourhood of $v'$ in
$\partial\mathbb H_\mathbb R^2$. In other words,
$A\subset\partial\mathbb H_\mathbb R^2$ is an open arc containing $v'$
and not containing $v$. By \ref{fixed_points_density_Lemma}, for a
suitable $p\in A\cap F$ and for some $1\ne\gamma\in\Gamma$, we have
$\gamma(p)=p$.

\begin{wrapfigure}[11]{r}{0.41\textwidth}
\includegraphics[clip=true,trim=0 0 0
0.5,width=0.41\textwidth]{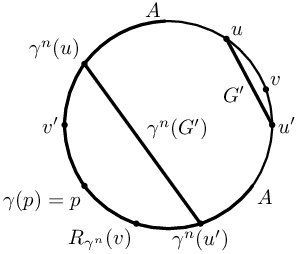}
\end{wrapfigure}
We consider two cases. The first case: $u,u'$ are not fixed by
$\gamma$. Then, taking into account that $\gamma$ is hyperbolic or
parabolic, we conclude by \ref{limit_Remark} that $\gamma^n(u)\to p$
and $\gamma^n(u')\to p$ for $n\to\infty$ or for $n\to-\infty$. Hence,
$\gamma^n(u),\gamma^n(u')\in A$ for some $n\in\mathbb Z$. Therefore,
$R_{\gamma^n}(v)\in A$.

The second case: one of $u,u'$ is fixed by $\gamma$. By
\ref{no_fix_point_Lemma}, there exists $\gamma_0\in\Gamma$ such that
the points $\gamma_0(u),\gamma_0(u')$ are not fixed by $\gamma$. Now,
by \ref{limit_Remark}, we have $\gamma^n\gamma_0(u)\to p$ and
$\gamma^n\gamma_0(u')\to p$ for $n\to\infty$ or for $n\to-\infty$. This
implies $R_{\gamma^n\gamma_0}(v)\in A$ for some $n\in\mathbb Z$.

For an arbitrarily small open arc $A$ containing $v'$, we found, in
either case, some $\gamma'\in\Gamma$ such that $R_{\gamma'}(v)\in A$
and $v\notin\gamma'(G')$. This implies \ref{hyperbolic_Lemma}.
\endproof

\hfill

\remark\label{negative_Remark}
{\sl We stated \ref{_lattices_dense_Proposition_}
in assumption that $(l,l)> 0$ (this assumption was geometrically
motivated). But, in fact, this assumption
is completely unnecessary. Moreover, as the
following result implies, 
the proof of \ref{_lattices_dense_Proposition_} becomes much
easier when\/ $(l,l)\leq 0$.}

\hfill

%%%%%%%%%%%%%%%%%%%%%%%%%%%%%%%%%%%%%%%%%%%%%%%%
\proposition\label{_negative_Proposition_}
{\sl The condition\/ $\langle l,l\rangle>0$ in
\ref{_lattices_dense_Proposition_} is unnecessary.}

\hfill

{\bf Proof:} To see this, we repeat the proof of
\ref{_lattices_dense_Proposition_} literally until
\ref{_3-dim-reformula_Lemma_}. To obtain \ref{negative_Remark}, we need
to check a version of \ref{_3-dim-reformula_Lemma_} when the vector $l$
is not assumed to be positive.

We reduce the case of $\langle l,l\rangle<0$ to the case of
$\langle l,l\rangle=0$. Let $l_0$ be a limit point of the orbit
$\Gamma\cdot l$. Since $\Gamma$ is a discrete subgroup in
$\PSL_2(\mathbb R)$, this limit lies on the absolute, and we have
$\langle l_0,l_0\rangle=0$. It~suffices to show that any geodesic $G$
passing through $l_0$ lies in the closure of the set of all geodesics
$G'$ that pass through $\gamma(l)$ for some $\gamma\in\Gamma$. For a
given point $\gamma(l)$, we denote by $G'_\gamma$ the Euclidean
parallel to $G$ passing through $\gamma(l)$. For this choice of
$G'_\gamma$, the  limit $\gamma(l)\to l_0$ implies the limit
$G'_\gamma\to G$.

It remains now to prove \ref{_3-dim-reformula_Lemma_} when
$\langle l,l\rangle=0$. Since $F$ is dense in
$\partial\mathbb H_\mathbb R^2$ by \ref{fixed_points_density_Lemma},
the subset $\Gamma\cdot l$ is also dense in
$\partial\mathbb H_\mathbb R^2$. So, fixing one end of an arbitrary
geodesic $G\in\Gr_{+-}(V)$, we can approximate the other one by a point
in $\Gamma\cdot l$.
\endproof

\hfill

{\bf Acknowledgements:}
We are grateful to Eyal Markman and Vyacheslav Nikulin for an email
correspondence.

{\small

\noindent
{\sc Misha Verbitsky}\\
{\sc  Laboratory of Algebraic Geometry, NRU-HSE,\\
7 Vavilova Str. Moscow, Russia, 117312}\\
{\tt verbit@maths.gla.ac.uk, \ \  verbit@mccme.ru\\}

\vskip5pt

\noindent
{\sc Sasha Anan$'$in}\\
{\sc IMECC -- UNICAMP, Departamento de Matem\'atica,\\
Caixa Postal 6065 13083-970 Campinas-SP, Brazil}\\
{\tt ananin\_sasha@yahoo.com}}


\begin{thebibliography}{VGSh}

\bibitem[AGr]{_AG:Geometries_}
Anan$'$in, S., Grossi, C.~H.,
{\em Coordinate-free classic geometries,}
Mosc.~Math.~J.~{\bf11} (2011), 633--655, arXiv:math/0702714

\bibitem[Bea]{_Beauville_}
Beauville, A.,
{\em Vari\'et\'es K\"ahleriennes dont la premi\`ere classe de Chern est
nulle,}
J.~Differential Geom.~{\bf18} (1983), no.~4, 755--782.

\bibitem[Bes]{_Besse:Einst_Manifo_}
Besse, A.~L.,
{\em Einstein Manifolds,}
Springer-Verlag, Berlin, Heidelberg, New York 1987, 516 pp.

\bibitem[Bo1]{_Bogomolov:decompo_}
Bogomolov, F.~A.,
{\em On the decomposition of K\"ahler manifolds with trivial canonical
class,}
Matt.~USSR-Sb.~{\bf22} (1974), no.~4, 580--583.

\bibitem[Bo2]{_Bogomolov:defo_}
Bogomolov, F.~A.,
{\em Hamiltonian K\"ahler manifolds,}
Sov.~Math.~Dokl. {\bf19} (1978), 1462--1465.

\bibitem[Bou]{_Boucksom_}
Boucksom, S.,
{\em Higher dimensional Zariski decompositions,}
Ann.~Sci.~Ecole Norm.~Sup. (4){\bf37} (2004), no.~1, 45--76,
arXiv:math/0204336

\bibitem[DP]{_Demailly_Paun_}
Demailly, J.-P., Paun, M.,
{\em Numerical characterization of the K\"ahler cone of a compact
K\"ahler manifold},
Ann. of Math.~{\bf159} (2004), 1247--1274,
arXiv:math/0105176

\bibitem[F]{_Fujiki:HK_}
Fujiki, A.,
{\em On the de Rham Cohomology Group of a Compact K\"ahler Symplectic
Manifold,}
Adv.~Stud.~Pure Math.~{\bf10} (1987), 105--165.

\bibitem[GHS1]{_GHK:moduli_HK_}
Gritsenko, V., Hulek, K., Sankaran, G.~K.,
{\em Moduli spaces of irreducible symplectic manifolds,}
Compositio Mathematica {\bf146} (2010), no.~2, 404--434,
arXiv:0802.2078

\bibitem[GHS2]{_GHK:modular_}
Gritsenko, V., Hulek, K., Sankaran, G.~K.,
{\em Abelianisation of orthogonal groups and the fundamental group of
modulae varieties,}
arXiv:0810.1614, 21 pp.

\bibitem[H1]{_Huybrechts:basic_}
Huybrechts, D.,
{\em Compact hyper-K\"ahler manifolds\/{\rm:} basic results,}
Invent.~Math.~{\bf135} (1999), no.~1, 63--113,
arXiv:alg-geom/9705025

\bibitem[H2]{_Huybrechts:finiteness_}
Huybrechts, D., 
{\em Finiteness results for hyperk\"ahler manifolds},
J.~Reine Angew. Math. {\bf558} (2003), 15--22,
arXiv:math/0109024



\bibitem[KS]{_Kontsevich-Soibelman:torus_}
Maxim Kontsevich, Yan Soibelman,
{\em Homological mirror symmetry and torus fibrations},
arXiv:math/0011041, Symplectic geometry and mirror
symmetry (Seoul, 2000), World Sci. Publishing, 
River Edge, NJ, 2001, pp. 203-263.

\bibitem[KV]{_Kamenova_Verbitsky_}
Ljudmila Kamenova, Misha Verbitsky,
{\em Families of Lagrangian fibrations on hyperkaehler manifolds},
arXiv:1208.4626, 13 pages.

\bibitem[M]{_Markman:constra_}
Markman, E.,
{\em Integral constraints on the monodromy group of the hyperk\"ahler
resolution of a symmetric product of a\/ {\rm K3} surface,}
International Journal of Mathematics {\bf21} (2010), no.~2, 169--223,
arXiv:math/0601304

\bibitem[T]{_Todorov:torsion_}
Todorov, A.,
{\em Ray Singer Analytic Torsion of CY Manifolds\/ {\rm II},}
arXiv: math/0004046, 18 pp.

\bibitem[V0]{_Verb:alge_} 
Verbitsky, M.,
{\it Algebraic structures on hyperk\"ahler manifolds,}
Math. Res. Lett., {\bf3} (1996), 763--767.

\bibitem[V1]{_Verbitsky:Symplectic_II_}
Verbitsky, M., {\em Trianalytic subvarieties 
of hyperkaehler manifolds}, GAFA vol. 5 no. 1 (1995) pp. 92-104.


\bibitem[V2]{_V:Torelli_}
Verbitsky, M.,
{\em A global Torelli theorem for hyperk\"ahler manifolds,}
arXiv: 0908.4121, 47 pp.


\bibitem[V3]{_Verbitsky:SYZ_}
Verbitsky, M., 
{\em Hyperkahler SYZ conjecture and semipositive line bundles},
 arXiv:0811.0639, GAFA 19, No. 5, 1481-1493 (2010)

\bibitem[V4]{_Verbitsky:parabolic_}
Verbitsky, M., 
{\em  Parabolic nef currents on hyperkaehler manifolds},
arXiv:0907.4217, 19 pages.

\bibitem[VGSh]{_VGSh:VINITI_}
Vinberg, E.~B., Gorbatsevich, V.~V., Shvartsman, O.~V.,
{\em Discrete Subgroups of Lie Groups,}
Encyclopaedia of Mathematical Sciences {\bf21}, Springer-Verlag, 2000,
224 pp.

\bibitem[Vi]{_Viehweg:moduli_}
Viehweg, E.,
{\em Quasi-projective Moduli for Polarized Manifolds,}
Springer-Verlag, Berlin, Heidelberg, New York, 1995,
Ergebnisse der Mathematik und ihrer Grenzgebiete, 3.~Folge, Band 30,
also available at
{\tt http://www.uni-due.de/$\widetilde{\phantom{a}}$mat903/books.html}

\end{thebibliography}
\end{document}